\documentclass[11pt]{article}
\usepackage{amssymb}
\usepackage{amsmath}
\usepackage{graphicx}
\usepackage{setspace}
\begin{document}

\newtheorem{theorem}{Theorem}[section]
\newtheorem{tha}{Theorem}
\newtheorem{conjecture}[theorem]{Conjecture}
\newtheorem{corollary}[theorem]{Corollary}
\newtheorem{lemma}[theorem]{Lemma}
\newtheorem{claim}[theorem]{Claim}
\newtheorem{proposition}[theorem]{Proposition}
\newtheorem{construction}[theorem]{Construction}
\newtheorem{definition}[theorem]{Definition}
\newtheorem{question}[theorem]{Question}
\newtheorem{problem}[theorem]{Problem}
\newtheorem{remark}[theorem]{Remark}
\newtheorem{observation}[theorem]{Observation}

\newcommand{\ex}{{\mathrm{ex}}}

\def\endproofbox{\hskip 1.3em\hfill\rule{6pt}{6pt}}
\newenvironment{proof}%
{%
\noindent{\it Proof.}
}%
{%
 \quad\hfill\endproofbox\vspace*{2ex}
}
\def\qed{\hskip 1.3em\hfill\rule{6pt}{6pt}}
\def\ex{{\rm\bf ex}}
\def\cA{{\cal A}}
\def\cB{{\cal B}}
\def\cF{{\cal F}}
\def\cG{{\cal G}}
\def\cH{{\cal H}}
\def\ck{{\cal K}}
\def\cI{{\cal I}}
\def\cJ{{\cal J}}
\def\cL{{\cal L}}
\def\cM{{\cal M}}
\def\cP{{\cal P}}
\def\cQ{{\cal Q}}
\def\cS{{\cal S}}

\def\cE{{\mathcal E}}

\voffset=-0.5in

\pagestyle{myheadings}
\markright{{\small \sc Frankl and F\"uredi:}
  {\it\small Another proof of the EKR theorem}}

\title{\huge\bf A new short proof of the EKR theorem}

\author{Peter Frankl\thanks{Peter Frankl Office, Shibuya-ku, Shibuya 3--12--25, Tokyo, Japan.
E-mail: peter.frankl@gmail.com} \quad and \quad
Zolt\'an F\"uredi\thanks{Dept.~of Mathematics, University of Illinois,
Urbana, IL 61801, USA. E-mail: z-furedi@illinois.edu
\newline
Research supported in part by the Hungarian National Science Foundation
 OTKA,  by the National Science Foundation under grant NFS DMS 09-01276,
  and by the European Research Council Advanced Investigators Grant 267195.
\newline
This copy was printed on \today, %
   \quad    {\rm\small {\jobname}.tex,}
   \hfill
Version as of
    August 10, 2011.
}}

\date{}

\maketitle
\begin{abstract}
A family $\cF$ is \emph{intersecting} if $F\cap F'\neq \emptyset$ whenever $F, F'\in \cF$.
Erd\H os, Ko, and Rado~\cite{EKR} showed that
\begin{equation}\label{eq:1}
|\cF|\leq {n-1\choose k-1}
   \end{equation}
  holds for an intersecting family of $k$-subsets of $[n]:=\{
  1,2,3,\dots,n\}$, $n\geq 2k$.
For $n> 2k$ the only extremal family consists of all $k$-subsets
 containing a fixed element.
Here a new proof is presented.
It is even shorter than the
 classical proof of Katona~\cite{Kat_EKR} using cyclic
 permutations, or the one found by Daykin~\cite{Day} applying the Kruskal-Katona
 theorem.
\end{abstract}

\section{Definitions: shadows, $b$-intersecting families.}
  ${X \choose k}$ denotes the family of $k$-element subsets of $X$.
For a 
  family of sets $\cA$ its $s$-\emph{shadow} $\partial_s \cA$ denotes
  the family of $s$-subsets of its members
$
  \partial_s \cA:= \{ S: |S|=s, \, \exists A\in \cA, S\subseteq A \}.
  $ 
E.g., $\partial_1\cA=\cup \cA$.
Suppose that $\cA$ is a family of $a$-sets such that $|A\cap A'|\geq b\geq 0$ for all
 $A,A'\in \cA$. 
Katona~\cite{Kat} showed that then
\begin{equation}\label{eq:Kat}
   |\cA|\leq |\partial_{a-b}\cA|.
  \end{equation}
We show that this inequality quickly implies the EKR theorem.

\section{The Proof}
Let $\cF\subset {[n]\choose k}$ be intersecting.
Define a partition $\cF_0:=\{ F\in \cF: 1\notin F\}$,
 $\cF_1:=\{ F\in \cF: 1\in F\}$
 and define
 $\cG_1:=\{ F\setminus \{ 1\}: 1\in  F\in \cF\}$.
Consider $\cF_0$ as a family on $[2,n]$.
Its complementary family $\cG_0:=\{ [2,n]\setminus F: F\in \cF_0\}$ is
 $(n-1-k)$-uniform.
The intersection property of $\cF$ implies that $F_1\setminus \{ 1\}=G_1\in \cG_1$
 is not contained in any member of  $[2,n] \setminus F_0=G_0\in \cG_0$.
We obtain
$$
   \cG_1 \cap \partial_{k-1}\cG_0 = \emptyset.
 $$
The intersection size $|G\cap G'|$ of $G,G'\in \cG_0$ is at least $n-2k$
$$
|G\cap G'|=|\left([2,n]\setminus F\right)
             \cap \left([2,n]\setminus F'\right)|=
              (n-1)-2k+|F\cap F'|.
  $$
Then (\ref{eq:Kat}) gives (with $a=n-k-1$, $b=n-2k\geq 0$) that
 $|\cG_0|\leq |\partial_{k-1}\cG_0|$.
 Summarizing
\begin{equation}\label{eq:4}
 |\cF|=|\cF_1|+|\cF_0|=|\cG_1|+|\cG_0|
  \leq |\cG_1|+|\partial_{k-1}\cG_0|
   \leq {n-1\choose k-1}.  \hfill \qed
 \end{equation}

\noindent
{\bf Extremal families.}\enskip
Equality holds in (\ref{eq:Kat}) if and only if $a=b$, or $\cA=\emptyset$, or
 $\cA\equiv {[2a-b]\choose a}$.
Thus, for $n> 2k$, equality in (\ref{eq:4}) implies either $\cG_0=\emptyset$ and
 $1\in \cap \cF$, or $\cG_0\equiv {[2,n-1] \choose n-1-k}$ and $n\in \cap \cF$.

\section{Two algebraic reformulations}

Given two families of sets $\cA$ and $\cB$, the {\it inclusion matrix} $I(\cA, \cB)$ is a $0$-$1$
 matrix of dimension $|\cA|\times |\cB|$, its rows and columns are labeled by the members of
 $\cA$ and $\cB$, respectively, the element $I_{A,B}=1$ if and only if $A\supseteq B$.
In the case $\cF\subseteq 2^{[n]}$ the matrix $I(\cF, {[n]\choose 1})$ is the usual {\it
 incidence matrix} of $\cF$, and $I(\cF, {[n] \choose s})$ is the
 {\em generalized} incidence matrix of order $s$.

Suppose that $L$ is a set of non-negative integers, $|L|=s$,
 and for any two distinct members $A, A'$ of the family $\cA$ one has
 $|A\cap A'|\in L$. Such a family is called $L$-\emph{intersecting}.
The Frankl-Ray-Chaudhuri-Wilson~\cite{FW,RW} theorem states that 
 in the case of $\cA\subseteq {[n]\choose k}$, $s\leq k$
 the row vectors of the generalized incidence matrix
 $I(\cA, {[n] \choose s})$ are linearly independent.
Here the rows are taken as real vectors (in~\cite{RW}) or as vectors over
 certain finite fields (in~\cite{FW}).
Note that this statement generalizes (\ref{eq:Kat})
 with $L=[k-1]$, $s=k-1$.

\noindent
{\bf Matrices and the EKR theorem.}\enskip
Instead of using (\ref{eq:Kat}) one can prove directly
 (like in~\cite{FW,RW}) that the row vectors of
 the inclusion matrix
 $I\left( G_0\cup G_1, {[2,n]\choose k-1} \right)$
 are linearly independent.

 \noindent{\bf Linearly independent polynomials.}\enskip
One can define homogeneous, multilinear polynomials $p(F,\mathbf{x})$
 of rank $k-1$ with variables $x_2, \dots, x_n$
 \begin{equation*}
   p(F, \mathbf{x})=\begin{cases} \sum \{ x_S: S\subset [2,n]\setminus F, \, |S|=k-1\}
        & \text{for}\enskip 1\notin F\in \cF,\\
     x_{F\setminus \{ 1\}} & \text{for}\enskip 1\in F\in \cF, \end{cases}
   \end{equation*}
where $x_S:=\prod_{i\in S} x_i$.
To prove (\ref{eq:1}) one can show that these polynomials are linearly independent.

\section{Remarks}

The idea of considering the shadows of the complements (one of the main steps of
 Daykin's proof~\cite{Day}) first appeared in Katona~\cite{Kat} (page 334) in 1964.
He applied a more advanced version of his intersecting shadow theorem~(\ref{eq:Kat}),
  namely an estimate using $\partial_{a-b+1}\cA$.

Linear algebraic proofs are common in combinatorics, see the book~\cite{BF}.
For recent successes of the method concerning intersecting families
  see Dinur and Friedgut~\cite{DF1,DF2}.
There is a relatively short proof of the EKR theorem in~\cite{FHW}
 using linearly independent polynomials.
In fact, 
 our proof here can be considered
 as a greatly simplified version of that one.

Since the algebraic methods are frequently insensitive to the structure
 of the hypergraphs in question it is much easier to give an upper bound
 ${n \choose k-1}$  
 which 
 holds for {\it all}  $n$ and $k$.
To decrease this formula to ${n-1\choose k-1}$ requires further insight.
Our methods resembles those of  Parekh~\cite{OP} and  Snevily~\cite{Sne}
 who succeeded to handle this for various related
 intersection problems.

Generalized incidence matrices proved to be extremely useful, see, e.g.,
 the ingenious proof of Wilson~\cite{W} for another Frankl-Wilson
 theorem, namely the exact form of the
 classical Erd\H os-Ko-Rado theorem concerning the maximum size of a
 $k$-uniform, $t$-intersecting family on $n$ vertices.
They proved~\cite{frankl1, W} that the maximum size is exactly
 ${n-t \choose k-t}$ if and only if $n\geq  (t+1)(k-t+1)$.


\end{document}